# RBF-FD Method for Some Dispersive Wave Equations and Their Eventual Periodicity

**Marjan Uddin[1,*], Hameed Ullah Jan[1,*] and Muhammad Usman[2]**

**Abstract:** In this paper, we approximate the solution and also discuss the periodic behavior termed as eventual periodicity of solutions of (IBVPs) for some dispersive wave equations on a bounded domain corresponding to periodic forcing. The constructed numerical scheme is based on radial kernels and local in nature like finite difference method. The temporal variable is executed through RK4 scheme. Due to the local nature and sparse differentiation matrices our numerical scheme efficiently recovers the solution. The results achieved are validated and examined with other methods accessible in the literature.

## 1 Introduction

The applications of dispersive wave equations are found in numerous branches of physical sciences from fluid dynamics, quantum mechanics, plasma physics to non-linear optics and to this extent in chemistry and biology as well [Linares and Ponce (2014); Pava (2009)]. As an example, the BBM equation developed by Benjamin, Bona and Mahony in the year 1972, also known as the regularized long-wave equation (RLWE) as a model illustrating the unidirectional propagation of long waves with small amplitude. The BBM model is well known in physical applications, see for example the references [Benjamin, Bona and Mahony (1972); Goldstein and Wichnoski (1980); Zhang, Wei and Gao (2002); Singh, Gupta and Kumar (2011); Peregrine (1966)]. This equation is an alternate model for the KdV equation derived by Korteweg and de Vries in 1895 with the presumption of small wave-amplitude and large wave length. The KdV equation and its various modifications serve as the modeling equations in several physical problems, see for example [Pava (2009); Sulem and Sulem (1999); Newell and Moloney (1992); Ablowitz, Kruskal and Ladik (1979); Brenner and von Wahl (1981); Tao (2006); Uddin, Jan, Ali et al. (2016)].

[1] Department of Basics Sciences, UET, Peshawar, Pakistan.
[2] Department of Mathematics, University of Dayton, Dayton, USA.
* Corresponding Authors: Marjan Uddin. Email: marjan@uetpeshawar.edu.pk;
Hameed Ullah Jan. Email: huj@uetpeshawar.edu.pk.

In the last century, in applied mathematics, physics and other related disciplines a rich area of research behind the KdV equation and its solution approximation has been developed see for example the references [Bona and Smith (1975); Bona, Pritchard and Scott (1981); Bona and Winther (1989); Fornberg and Whitham (1978); Kumar, Singh, Kumar et al. (2015); Al-Khaled (2001); Goswami, Singh and Kumar (2017); Savaşaneril and Hacioğlu (2018); Yüzbaşı and Şahin (2013); Alquran and Al-Khaled (2011)].

Another certain qualitative aspects of solutions of some dispersive wave equations indicated through experiments, which are linked together with their large-time behavior termed as eventual time periodicity of initial-boundary-value problems (IBVPs) solutions. A piston-type or flap-type wave maker fitted at the end of a channel in laboratory experiments show this attractive event. When the wavemaker periodically oscillates with a period $T$, then it is observed that the amplitude of the wave becomes periodic of period $T$ at each point along the channel after some time [Bona, Pritchard and Scott (1981)]. This interesting phenomenon of eventual periodicity has been elaborated previously in separate works, for the generalized BBM and KdV equations and their dissipative counterparts respectively which include Burger-type term investigated by Bona et al. [Bona and Wu (2009)]. The initial and boundary conditions will be appended with all these equations. In the recent work [Usman (2007); Shen, Wu and Yuan (2007)] a new solution representation and also the eventual periodicity is re-established for the KdV equation. The forced oscillations of the KdV equation and its stability and numerical treatment have been carried out very recently [Usman and Zhang (2009); Khan and Usman (2012); Al-Khaled, Haynes, Schiesser et al. (2018)]. In the present work we investigate this behavior using RBF-FD method for the BBM and KdV equations and of the equations acquired by take together two unlike dissipation with the aforementioned equations. More precisely, we estimate the solutions of the following IBVPs for the BBM and KdV equations

$$u_t + u_x + uu_x - u_{xxt} = 0, x, t \geq 0, \tag{1}$$

$$u_t + u_x + uu_x - u_{xxx} = 0, x, t \geq 0, \tag{2}$$

respectively, and also take their dissipative counterparts namely

$$u_t + u_x + uu_x - u_{xxt} - vu_{xx} = 0, x, t \geq 0, \tag{3}$$

$$u_t + u_x + uu_x - u_{xxx} - vu_{xx} = 0, x, t \geq 0, \tag{4}$$

where $v > 0$. The following conditions will be appended to all these equations

$$u(x, 0) = u_0(x), x \geq 0, u(0, t) = g(t), t \geq 0. \tag{5}$$

subject to the condition that the function $g$ will be a periodic of period $T > 0$, in other words $g(t + T) = g(t)$, where $t \geq 0$.

The radial basis functions (RBFs) approach is generally developed by Hardy [Hardy (1971)], the said work is carried out for interpolating 2D scattered data. The RBFs approach is the most utilized tool in the field of multivariate approximation theory.

The remarkable properties of RBFs include higher order smoothness, powerful convergence and efficiency in practical problems. Image processing, cartography, neural network, meteorology and turbulence analysis are some fields where RBFs methods have applications see for example the references [Buhmann (2000, 2003)]. A long list of mathematical applications of RBFs used in numerical methodologies for solving PDEs with high accuracy in multi-dimensions can be found in many fields of applied sciences see for example [Fasshauer (2007)]. These finding speedily developed research in RBFs and RBFs methods possessed considerable attention in scientific community as a truly mesh-free methods and of their ability to achieve spectral accuracy for the PDEs solutions on irregular domain contrast to other state-of-arts methods [Belytschko, Krongauz, Organ et al. (1996); Buhmann (2003)]. Some of the drawbacks and difficulties (like resulting linear system ill-conditioning) of RBFs methods have also been resolved through different techniques.

The RBF-FD technique is a combination of RBFs and conventional finite differences (FD) to get best accuracy as compare to PS and global RBFs techniques on scattered nodes without requiring a computational mesh. In particular, the RBF-FD method can offer high computational efficiency, because of the local approximations. This approach leads to sparse differentiation matrices rather than to full dense matrix problems, and can also be used for local (adaptive) refinements, and equally incorporating for large-scale computations in parallel. Originally this approach was developed by Tolstykh in the year 2000 [Tolstykh (2000)]. It has been modified and applied successfully to a variety of problems in recent years [Wang and Liu (2002); Wright (2003); Lee, Liu and Fan (2003); Chantasiriwan (2004); Davydov and Oanh (2011); Kosec and Zinterhof (2013); Shan, Shu and Qin (2009); Uddin, Ali and Ali (2015); Mramor, Vertnik and Šarler (2013)]. In this work, we construct a numerical scheme based on RBF-FD local method for the solution of the model Eqs. (1)-(5).

## 2 Description of the method

For mathematical formulation, we consider a general time dependent partial differential equation and describe gradually the RBF-FD method. Our model problem is of the type

$$\frac{\partial u(x,t)}{\partial t} = \mathfrak{L}u(x,t), \text{ such that } x \in \Omega \subset \mathbb{R}^d, d \geq 1, t > 0, \tag{6}$$

subject to the following initial condition and well define boundary conditions

$$u(x,0) = f(x), \mathfrak{B}u(x) = g(x,t), x \in \partial\Omega, \tag{7}$$

where the spatial operators $\mathfrak{L}$, $\mathfrak{B}$ denote the differential operators, and $f$, $g$ are given functions. Suppose $N$ denote the nodes in the domain $\Omega$ of the given problem used for approximation, then the RBF-FD method approximate $\mathfrak{L}u(x,t)$ at every single inner node as a linear combination at nodes having neighboring locations. To approximate $\mathfrak{L}u(x,t)$ at node say $x_1 x_2, x_3, \ldots, x_n$ are some $n$–1 nearest neighbor nodes of $x_1$, known as support region or stencil of n nodes for the node $x_1$ [Fasshauer and McCourt (2015)].

### 2.1 RBFs global differentiation matrix

The discretization of Eqs. (6)-(7) through RBF-FD method can be carried out by approximating $u$ using $\phi$ as a radial kernel at the node $x$ and is given by

$$\hat{u}(x) = \sum_{j=1}^{N} c_j \phi(\|x - x_j\|) = \mathbf{\Phi}(\mathbf{x})^T \mathbf{c}, x \in \Omega, \tag{8}$$

such that $\mathbf{\Phi}(\mathbf{x})^T = (\phi\|x - x_1\|, \phi\|x - x_2\|, \ldots, \phi\|x - x_N\|)$, and $\mathbf{c}$ is an unknown coefficient vector. The Lagrange form of Eq. (8) is given by

$$\hat{u}(x) = \mathbf{\Phi}(\mathbf{x})^T \mathbf{K}^{-1} u, \tag{9}$$

here $\mathbf{K}$ stand for the global RBFs system interpolation matrix. Since the kernel-based interpolant $\hat{u}$ in Eq. (9) yields a fine approximation to $u$. So any operator applied on $\hat{u}$ will also be a nice approximation of the corresponding operator applied to $u$ (see [Fasshauer and McCourt (2015); Fasshauer (2007))]. Application of a linear differential operator $\mathfrak{L}$ to the above Eq. (9) gives

$$\mathfrak{L}\hat{u}(x) = \mathfrak{L}\mathbf{\Phi}(\mathbf{x})^T \mathbf{K}^{-1} u. \tag{10}$$

From Eq. (10) we used the following notation for the values

$$\mathbf{K}_{\mathfrak{L}} = \begin{bmatrix} \mathfrak{L}\mathbf{\Phi}(\mathbf{x_1})^T \\ . \\ . \\ . \\ \mathfrak{L}\mathbf{\Phi}(\mathbf{x_N})^T \end{bmatrix}. \tag{11}$$

Thus the global differentiation (discretization) matrix $\mathbf{L}$ of order $N\times N$ matrix can be seen as

$$\mathbf{L} = \mathbf{K}_{\mathfrak{L}} \mathbf{K}^{-1}. \tag{12}$$

Since from Eq. (12), we see that the $i^{th}$ row of $\mathbf{K}_{\mathfrak{L}}$ corresponds to $\mathfrak{L}\mathbf{\Phi}(\mathbf{x_i})^T$ therefore we observe from Eq. (12) that the $i^{th}$ row of $\mathbf{L}$,

$$\mathbf{L}_i = \mathfrak{L}\mathbf{\Phi}(\mathbf{x_i})^T \mathbf{K}^{-1}, \tag{13}$$

represents one row of the global differentiation matrix $\mathfrak{L}$.

### 2.2 RBFs local differentiation matrix

The weights of local finite differences with respect to the point $x_i$, produce a local interpolant in a small local neighborhood of $x_i$ called a stencil denoted as $N_{xi}$. Let us define a set of points $\Omega = \{x_1, \ldots, x_N\}$ at which we want the approximation of derivatives, these points may be consider as stencil centers. Let $x_i$ be the $i^{th}$ evaluation node, then $n$ represent the number of nearest neighbor nodes of $x_i$ in the stencil $N_{xi}$, also define the set of points where we want to sample data $Z = \{z_1, \ldots, z_N\}$, the points in the stencil of size $n$ are collected in $Z_i \subset Z$. Now the approximation for the differential operator $L$ at the stencil having center at $x_i$ and collected at $Z_i \subset Z$ is defined as

$$\mathbf{L}_i = \mathbf{K}_{\mathfrak{L}}^{x_i}\mathbf{K}_{Z_i}^{-1}. \tag{14}$$

In fact it concerns only to the local stencil around $x_i$, therefore we name it a local differentiation matrix but its acts globally as it utilize the whole data of that stencil collectively, these $\mathbf{L}_i$ matrices have all non-zeros entries in the global sparse matrix $\mathbf{L}^{FD}$, but further it need to fix their locations in that matrix.

The $i^{th}$ row of $\mathbf{L}_i^{FD}$ contains non-zero entries in the matrix $\mathbf{L}_i$ this is actually a row vector as it has one evaluation node $x_i$. The nodes in $Z_i$ can be used to construct $\mathbf{L}_i$. To arrange correctly those points in the sparse row $\mathbf{L}_i^{FD}$ we can define an incidence matrix $\mathbf{P}_i \in \{0,1\}^{N_{x_i} \times N}$. Using this, the full sparse matrix $\mathbf{L}^{FD}$ is given by

$$\mathbf{L}^{FD} = \begin{bmatrix} \mathbf{K}_{\mathfrak{L}}^{x_1}\mathbf{K}_{Z_1}^{-1}\mathbf{P}_1 \\ \cdot \\ \cdot \\ \cdot \\ \mathbf{K}_{\mathfrak{L}}^{x_N}\mathbf{K}_{Z_N}^{-1}\mathbf{P}_N \end{bmatrix}. \tag{15}$$

where the entries of $\mathbf{P}_i$ are defined as

$$[\mathbf{P}_i]_{k,\ell} = \begin{cases} 1, & \text{if } k = \ell, \text{i.e., } k^{th} \text{ point in } Z_i \text{ matches the } \ell^{th} \text{ point in Z,} \\ 0, & \text{else.} \end{cases}$$

Hence it is clear that the $i^{th}$ row accepts values on $Z_i$. Finally the discretization of problem (6)-(7) is stated as

$$\acute{u} = \mathbf{M}u, \tag{16}$$

where $\mathbf{M} = \begin{bmatrix} \mathbf{L}^{FD} \\ \mathbf{B}^{FD} \end{bmatrix}$, in which $\mathbf{B}^{FD}$ rep resent the discretization of operator applied at the boundary and can be found accordingly as $\mathbf{L}^{FD}$. To evolve the ODE system (16) in time any ODE solver like ode23, ose45, ode113 and many other from Matlab can be used.

### *2.3 Stability of the proposed numerical scheme*

In our proposed numerical scheme which is based on RBF-FD method we transformed the time-dependent partial differential equation into an ODEs system in time. This type of technique is called the method of lines by which we can solve this system of coupled ODEs using the finite difference method in time for example Runge-Kutta methods, etc. The method of lines stability may be estimated by the well known rule of thumb. It is shown in the work [Trefethen (2000)], that the method of lines will be stable, when the eigenvalues of spatial discretization operator, linearized and scaled by step size $\delta t$, lie in region of stability of the corresponding time-discretization operator.

In the present RBF-FD method, the region of stability is a part of the complex plane containing those eigenvalues, where the technique produces a bounded solution. The RBF-FD meshless method of lines numerical scheme is defined in Eq. (16). We used the above criteria of satiability for our numerical scheme. It is shown that the current RBF-FD (localized) numerical scheme is unconditionally stable for all values of RBFs shape parameter and small step size $\delta t$, when solving the proposed model equations.

## 3 Application of the proposed numerical scheme for dispersive wave equations

First we test the RBF-FD numerical scheme to the KdV equation with known exact solution, and then apply to the model equations given in (1)-(5).

$$u_t + \beta u u_x + \gamma u_{xxx} = 0, x \in [a, b], t > 0, \tag{17}$$

$$u(x,t) = \frac{c \sec h^2 (\frac{\sqrt{c}}{2}(x - ct) - 7)}{2}. \tag{18}$$

We used the initial solution $u_0$ as well as the boundary conditions from the analytic solution given in Eq. (18). The problem in Eq. (17) is solved over the spatial domain [0,40], where as the time domain [0,5] and the parameters $\beta$=6, $c$=0.50 and $\gamma$=1. The accuracy and efficiency of the present scheme is verified in terms of the $L_\infty$ and $L_2$ error norms defined by

$$L_\infty = \|u^{ex} - u^{ap}\|_\infty = \max|u^{ex} - u^{ap}| \tag{19}$$

$$L_2 = \|u^{ex} - u^{ap}\|_2 = \sqrt{h \sum_{i=1}^{N} (u^{ex} - u^{ap})^2} \tag{20}$$

The results of RBF-FD method are compared with methods in literature and it is observed that the applied RBF-FD method is very beneficial with respect to computational time and also in terms of accuracy and convergence rate as evident from Tab. 1 and Fig. 1 respectively. The numerical stability of the proposed numerical scheme for problem (17) in the graphical form is shown in Fig. 1. It is evident that all the eigenvalues of the differentiation matrix are well kept in stability region of RK4 method. We also observed that the present numerical scheme is unconditionally stable for a large range of RBFs shape parameters and small step size $\delta t$. This is mainly due to the local and sparse nature of small size system matrices defined over the small local sub-domains called stencils.

**Table 1:** Comparison of RBF-FD method with other methods, when $c$=0.50, $\delta t$=0.001, $\beta$=6, $\gamma$=1, $C_{MQ}$=1.2, $x\in[0, 40]$, $N$=200, $N_x$=25, $t\in[0, 5]$, corresponding to (17)

| Method | t | $L_\infty$ | $L_2$ | C. time |
|---|---|---|---|---|
| [Uddin, Shah and Ali (2015)] | 1 | 1.2800e-006 | 5.2490e-006 | 5.320 |
| | 2 | 1.9340e-006 | 6.5130e-006 | 9.860 |
| | 3 | 2.3540e-006 | 7.9050e-006 | 14.570 |
| | 4 | 2.9030e-006 | 9.9260e-006 | 19.080 |
| | 5 | 3.8880e-006 | 1.3020e-005 | 24.280 |
| [Shen (2009)] | 1 | 1.804e-005 | 6.236e-005 | 14.00 |
| | 2 | 3.037e-005 | 1.126e-004 | 20.00 |
| | 3 | 4.008e-005 | 1.553e-004 | 25.00 |
| | 4 | 4.834e-005 | 1.940e-004 | 30.00 |
| | 5 | 5.609e-005 | 2.294e-004 | 36.00 |
| RBF-FD(MQ) method | 1 | 8.468e-005 | 4.803e-004 | 0.33 |
| | 2 | 1.538e-004 | 9.336e-004 | 0.61 |
| | 3 | 2.207e-004 | 1.347e-003 | 0.85 |
| | 4 | 2.692e-004 | 1.712e-003 | 1.10 |
| | 5 | 3.119e-004 | 2.032e-003 | 1.36 |

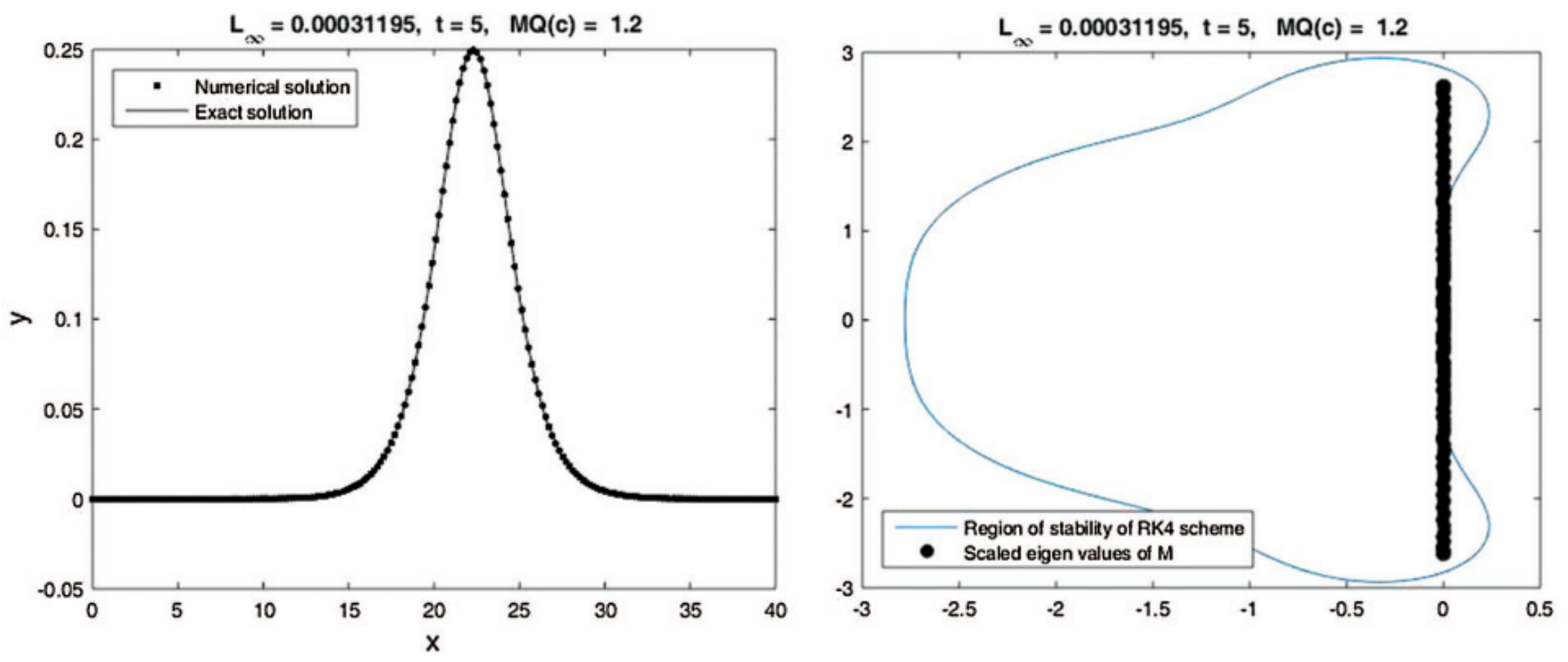

**Figure 1:** Solitary wave: Stable eigenvalue spectrum, $c$=0.50, $\delta t$=0.001, $\beta$=6, $\gamma$=1, $C_{MQ}$=1.2, $x\in[0, 40]$, $N$=200, $N_x$=25, $t\in[0, 5]$

## 4 Eventual periodicity of some dispersive wave equations

In this section we applied the RBF-FD method to study numerically the eventual periodicity of the BBM and KdV equations. We also study some other modified form by assemble two different kind of dissipation in these equations (see (1)-(5)). For a suitable boundary data $g$, we established the eventual periodicity in which the initial data $u_0$ is not absolutely necessary, hence we set it to zero simply. The eventual periodicity is studied numerically using the RBF-FD method solutions in graphical form. Since any solution $u(x,t)$ of these IBVPs converges to 0 as $x \to \infty$. We estimate the solutions of the following scaled problems studied by Shen et al. in [Shen, Wu and Yuan (2007)].

$$\begin{cases} u_t + \alpha u_x + \beta u u_x - \mu u_{xxt} = 0, & x \in [-1,1], \quad t \in [0,T], \\ u(0,t) = g(t), & \\ u(1,t) = u_x(1,t) = 0, & t \in [0,T], \\ u(x,0) = u_0(x), & x \in [-1,1], \end{cases} \tag{21}$$

and

$$\begin{cases} u_t + \alpha u_x + \beta u u_x - \mu u_{xxt} - \delta u_{xx} = 0, & x \in [-1,1], \quad t \in [0,T], \\ u(0,t) = g(t), & \\ u(1,t) = u_x(1,t) = 0, & t \in [0,T], \\ u(x,0) = u_0(x), & x \in [-1,1], \end{cases} \tag{22}$$

similarly

$$\begin{cases} u_t + \alpha u_x + \beta u u_x - \mu u_{xxx} = 0, & x \in [-1,1], \quad t \in [0,T], \\ u(0,t) = g(t), & \\ u(1,t) = u_x(1,t) = 0, & t \in [0,T], \\ u(x,0) = u_0(x), & x \in [-1,1], \end{cases} \tag{23}$$

and

$$\begin{cases} u_t + \alpha u_x + \beta u u_x - \mu u_{xxx} - \delta u_{xx} = 0, & x \in [-1,1], \quad t \in [0,T], \\ u(0,t) = g(t), & \\ u(1,t) = u_x(1,t) = 0, & t \in [0,T], \\ u(x,0) = u_0(x), & x \in [-1,1], \end{cases} \tag{24}$$

and also

$$\begin{cases} u_t + \gamma u + \alpha u_x + \beta u u_x - \mu u_{xxx} = 0, & x \in [-1,1], \quad t \in [0,T], \\ u(0,t) = g(t), & \\ u(1,t) = u_x(1,t) = 0, & t \in [0,T], \\ u(x,0) = u_0(x), & x \in [-1,1], \end{cases} \tag{25}$$

where $\gamma$, $\alpha$, $\beta$, $\mu$ and $\delta$ are parameters and $g(t)=\sin(20\pi t)\tanh(5t)$. For $u_0 \equiv 0$, Eqs. (21)-(25) are valid approximations of Eqs. (1)-(5) respectively, till the instant when the wave-front generated by the boundary data $g$ (periodic function of period $T$>0) arrives the right boundary point $x$=1.

### *4.1 Eventual periodicity for Linearized BBM equation*

The RBF-FD method is used for the numerical solutions of model Eq. (21) over the spatial domain [-1,1], and in a time domain [0,1.8]. In all of our computations we use multiquadric radial basis function defined by $\phi(r, C) = \sqrt{C^2 + r^2}$, contains a shape parameter $C$. The solution accuracy mainly depends on this shape parameter $C$. But in this particular local setting it does not much effect the solution as the localized RBF-FD method is stable for a large range of shape parameters $C$. For the linearized BBM equation a first set of our experiments with parameters $\alpha$=1.0, $\beta$=0 and $\mu=10^{-6}$. The amplitudes $u(x,t)$ recorded in the six graphs shown in Fig. 2 at particular points $x$=-0.950670, -0.808460, -0.587280, -0.308720, 0 and 0.999650 in the domain [-1,1]. Here $N$ denote total points in the domain, while $N_x$ denote the points in local sub-domain respectively. The horizontal and vertical axes stand for the time $t$ and amplitude $u$ respectively in these graphs. Initially the boundary data are not exactly periodic but after some time it becomes eventually periodic. In each problem the last graph at $x$=0.999650 shows the amplitude remain zero.

### *4.2 Eventual periodicity for Linear BBM-Burgers equation*

We compute the solutions of (22) by RBF-FD method, for a linear BBM-Burgers equation with parameters $\alpha$=1.0, $\beta$=0, $\mu=10^{-6}$ and $\delta=10^{-5}$. The amplitudes $u(x,t)$ of the model is recorded in the six graphs in Fig. 3 at particular points $x$=-0.950670, -0.808460, -0.587280, -0.308720, 0 and 0.999650 in the spatial domain [-1,1] for the time [0,1.8]. The plots in the figures clearly validate the eventual periodic behavior of the solution at these particular positions in the given domain. The effect of the Burgers type dissipation is observed in the damped amplitude.

### *4.3 Eventual periodicity for non-linear BBM equation*

Here in this problem the solution of the model Eq. (21) is obtained by RBF-FD method. In our experiment we used the parameters $\alpha$=1.0, $\beta$=0.05 and $\mu=10^{-6}$. We detect the form of eventual periodicity in $u(x,t)$ see Fig. 4 at all selected positions, and noticed the effect of the nonlinear term which increase the amplitude as compare to the plots with ones for the linearized BBM equation.

### *4.4 Eventual periodicity for nonlinear BBM-Burgers equation*

We solved the model Eq. (22) by RBF-FD method as a fourth set of our experiments with parameters $\alpha$=1.0, $\beta$=0.05, $\mu=10^{-6}$ and $\delta=10^{-5}$. Again we discovered the design of eventual periodicity in $u(x,t)$ at all selected positions, see Fig. 5, and the effect of Burgers type dissipation in the damped amplitude.

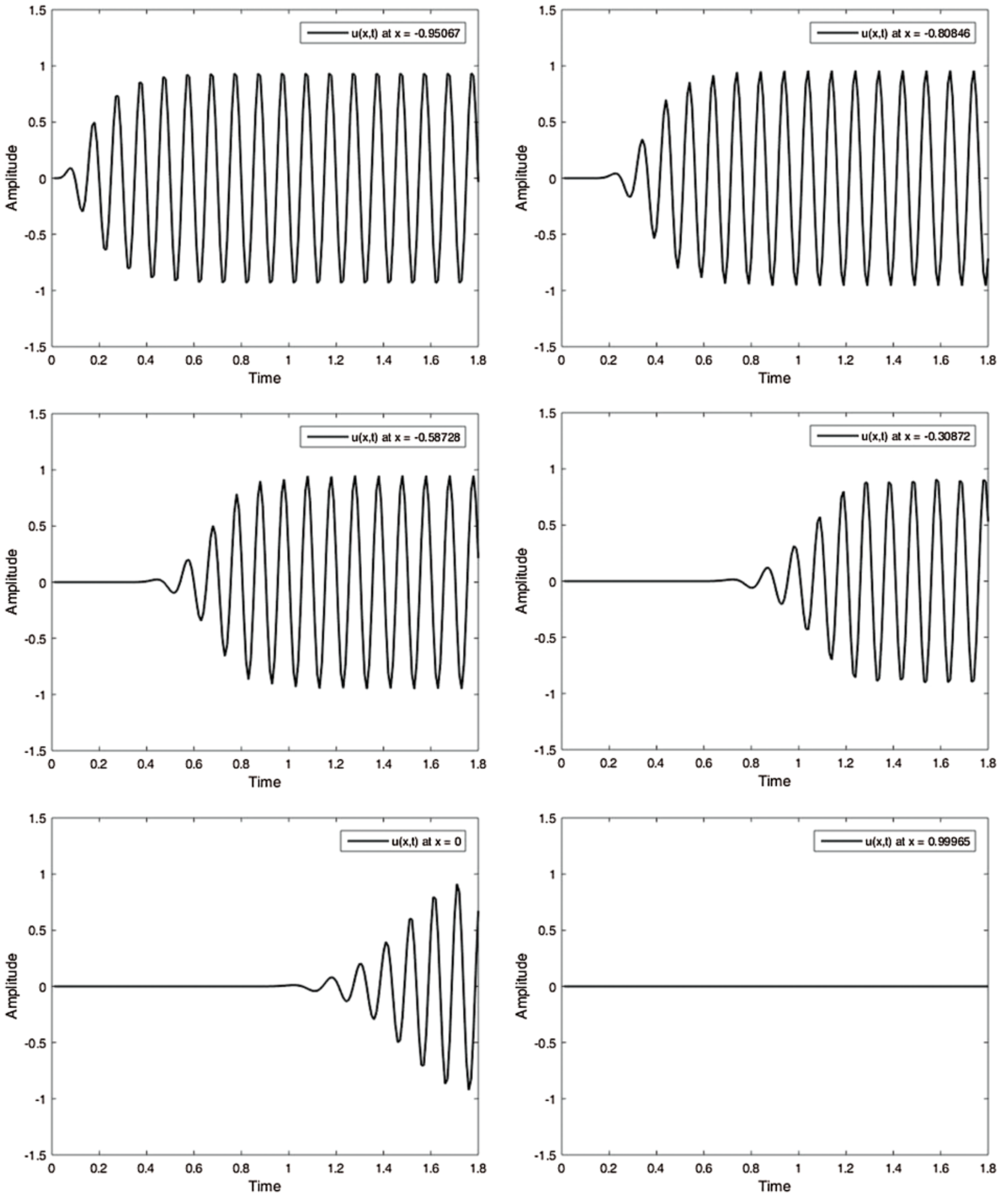

**Figure 2:** Eventual periodicity for linearized BBM equation: $u(x,t)$ at $x$=-0.950670, -0.808460, -0.587280, -0.308720, 0 and 0.999650 are shown through above six graphs respectively. For $\alpha$=1, $\beta$=0, $\mu=10^{-6}$, $x\in[-1, 1]$, $C_{MQ}$=0.001, $N$=200, $N_x$=25, $tmax$=1.8, $\delta t$=0.001, $g(t)=\sin(20\pi t)\tanh(5t)$

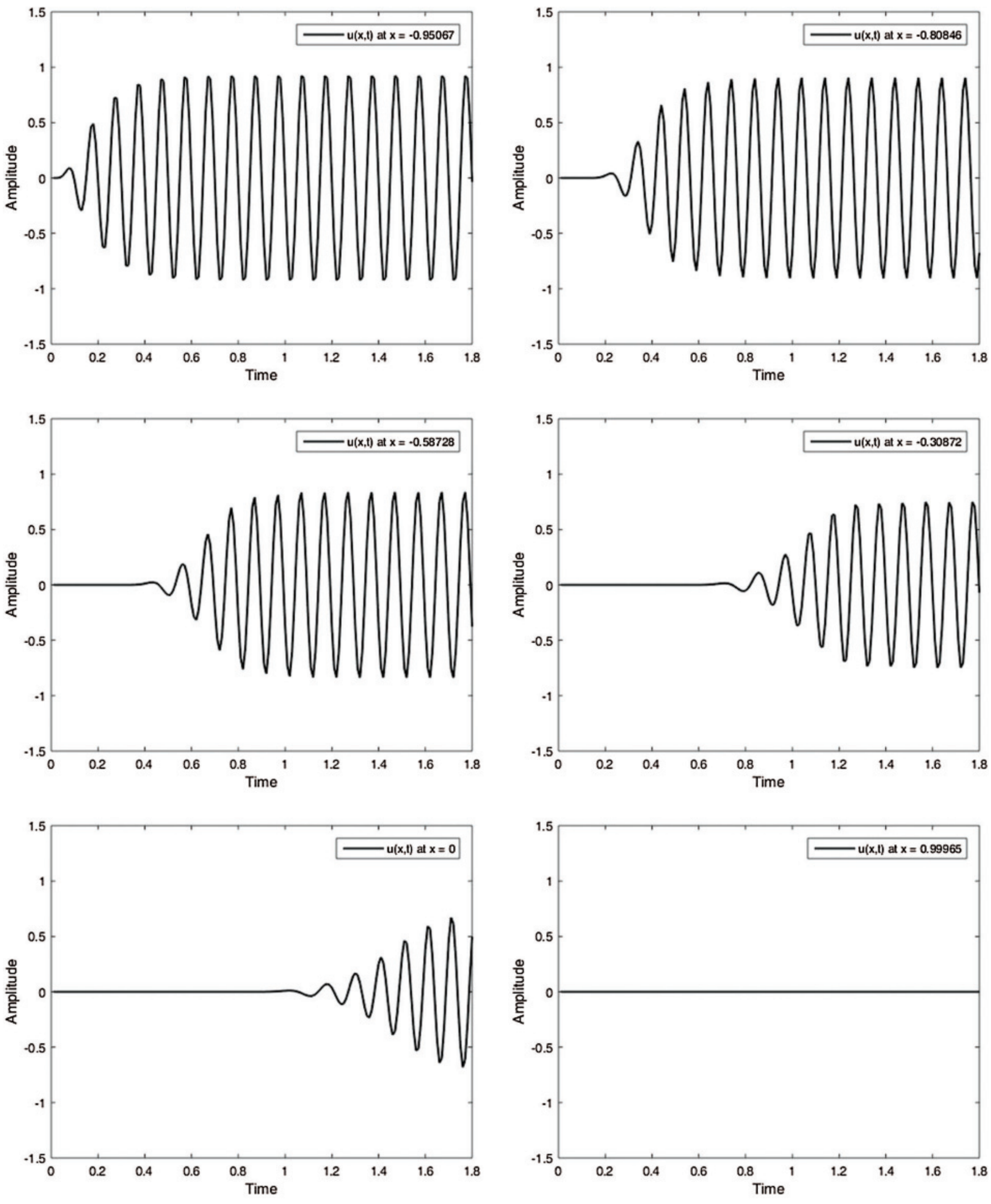

**Figure 3:** Eventual periodicity for linear BBM-Burgers equation: $u(x,t)$ at $x$=-0.950670, -0.808460, -0.587280, -0.308720, 0 and 0.999650 are shown through above six graphs respectively. For $\alpha$=1, $\beta$=0, $\mu$=$10^{-6}$, $\delta$=$10^{-5}$, $C_{MQ}$=0.001, $x \in [-1, 1]$, $N$=200, $N_x$=25, $tmax$=1.8, $\delta t$=0.001, $g(t)$=sin(20$\pi t$) tanh(5$t$)

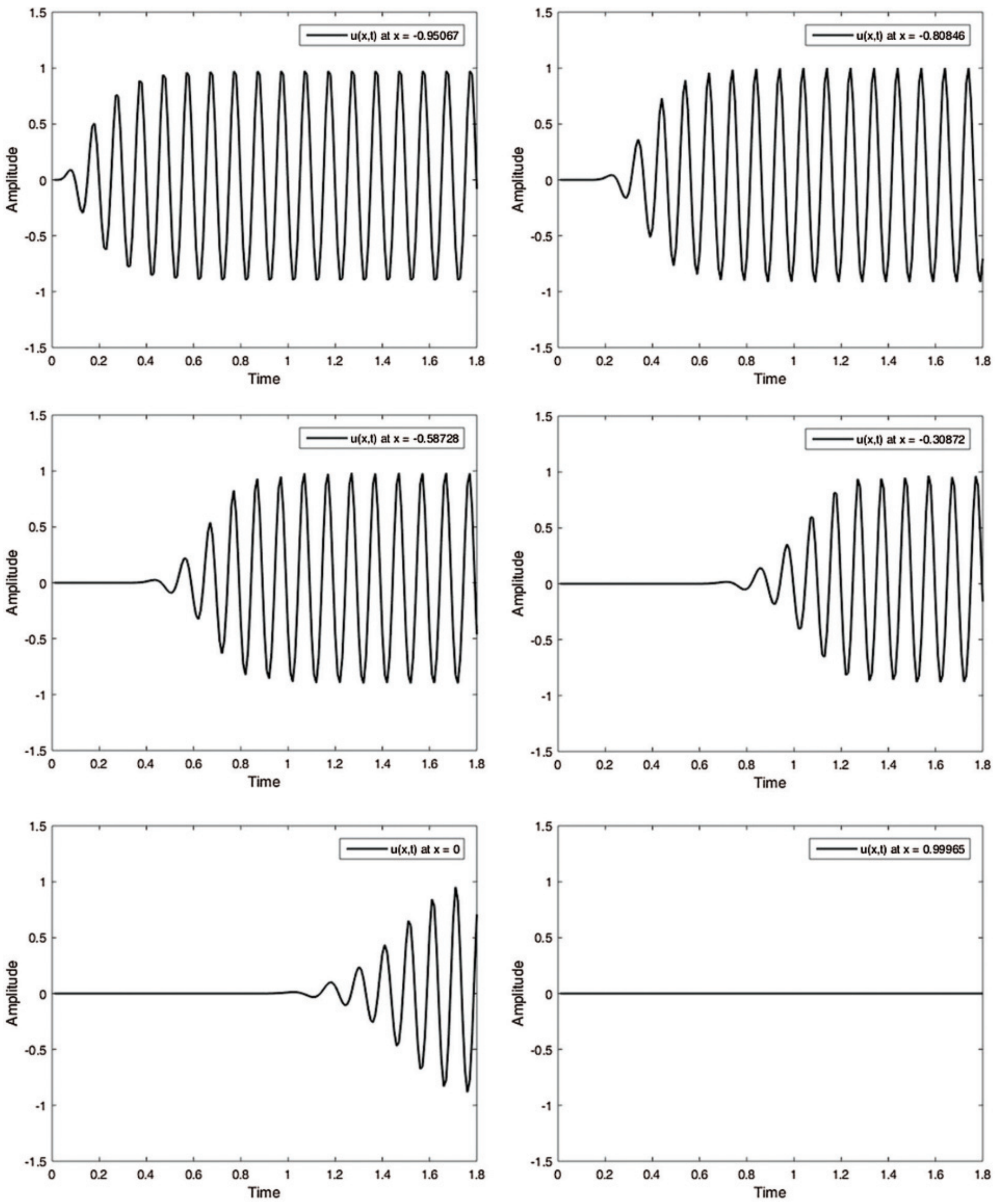

**Figure 4:** Eventual periodicity for non-linear BBM equation: $u(x,t)$ at $x$=-0.950670, -0.808460, -0.587280, -0.308720, 0 and 0.999650 are shown through above six graphs respectively. For $\alpha$=1, $\beta$=0.05, $\mu=10^{-6}$, $x\in[-1, 1]$, $C_{MQ}$=0.001, $N$=200, $N_x$=25, $tmax$=1.8, $\delta t$=0.001, $g(t)=\sin(20\pi t)\tanh(5t)$

## 5 Eventual periodicity for KdV equation

In this section, we obtained the solutions of model Eqs. (23)-(25), namely KdV type equations, using RBF-FD method. The six graphs in each experiments recorded the amplitudes $u(x,t)$ at the particular points $x$=-0.950670, -0.808460, -0.587280, -0.308720, 0 and 0.999650 in the given domain, and over the time interval $0<t<1.8$. The horizontal and vertical axes represents the time $t$ and amplitude $u$ respectively in these graphs. Initially the boundary data are not exactly periodic but after some time become eventually periodic. In each problem the last graph at $x$=0.999650 shows the amplitude remains zero.

### *5.1 Eventual periodicity for linearized KdV equation*

We consider the model Eq. (23) and solved it by RBF-FD method with parameters $\alpha$=1.0, $\beta$=0 and $\mu=10^{-5}$, Fig. 6 clearly validate the eventual periodicity of the solutions at all selected positions.

### *5.2 Eventual periodicity for linear KdV-Burgers equation*

We consider the model Eq. (24) with parameters $\alpha$=1.0, $\beta$=0, $\mu=10^{-5}$ and $\delta=10^{-4}$, and solved by RBF-FD method to study the eventual periodicity for the linearized KdV-Burgers equation as shown in Fig. 7 and observe the effect of Burgers type dissipation on the damped amplitude.

### *5.3 Eventual periodicity for non-Linear KdV equation*

Here the model Eq. (23) is solved by RBF-FD method with parameters $\alpha$=1.0, $\beta$=0.05 and $\mu=10^{-5}$, We computed the eventual periodicity for the non-Linear KdV equation, see for example Fig. 8. By comparing these plots with ones for the linearized KdV equation the effects of the nonlinear term as a significant increase in amplitude is observed.

### *5.4 Eventual periodicity for KdV-Burgers equation*

As fourth set of experiments, we consider the model Eq. (24) with parameters $\alpha$=1.0, $\beta$=0.05, $\mu=10^{-5}$ and $\delta=10^{-4}$, we compute the eventual periodicity for the KdV-Burgers equation, using RBF-FD method and the results are shown in Fig. 9. It is observed that the effect of Burgers type dissipation on the damped amplitude is evident.

Finally to realize the effects on the eventual periodicity of the damping term $\eta\, u$ for a general boundary data, we consider the model Eq. (25) with parameters $\eta$=4.5, $\alpha$=1.0, $\beta$=0.05 and $\mu=10^{-4}$. We noticed that the pattern of eventual periodicity is still remained but the amplitudes are significantly damped, see Fig. 10.

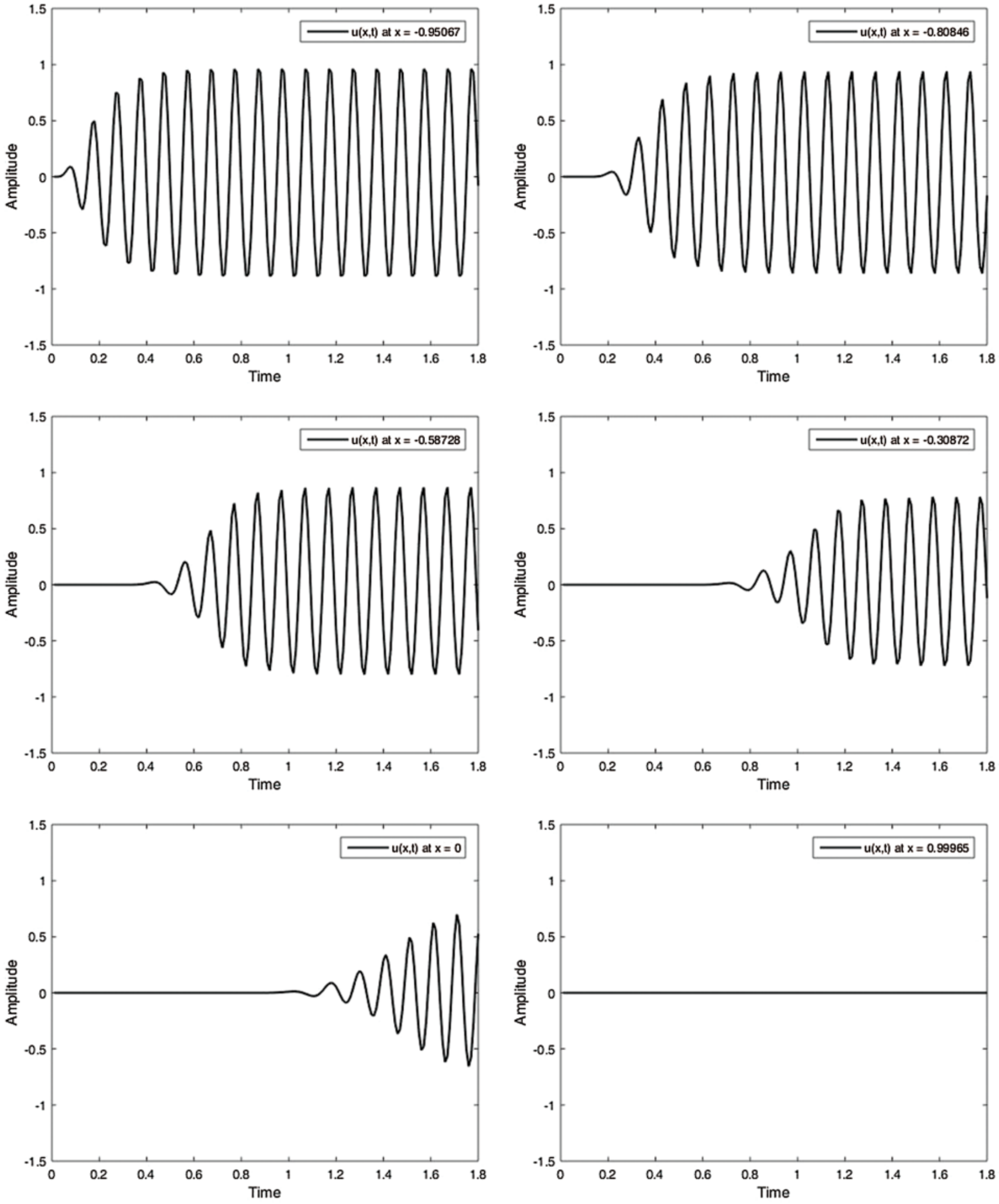

**Figure 5:** Eventual periodicity for nonlinear BBM-Burgers equation: $u(x,t)$ at $x$=-0.950670, -0.808460, -0.587280, -0.308720, 0 and 0.999650 are shown through above six graphs respectively. For $\alpha$=1, $\beta$=0.05, $\mu$=$10^{-6}$, $\delta$=$10^{-5}$, $x\in[-1, 1]$, $C_{MQ}$=0.001, $N$=200, $N_x$=25, $tmax$=1.8, $\delta t$=0.001, $g(t)=\sin(20\pi\ t)\ \tanh(5t)$

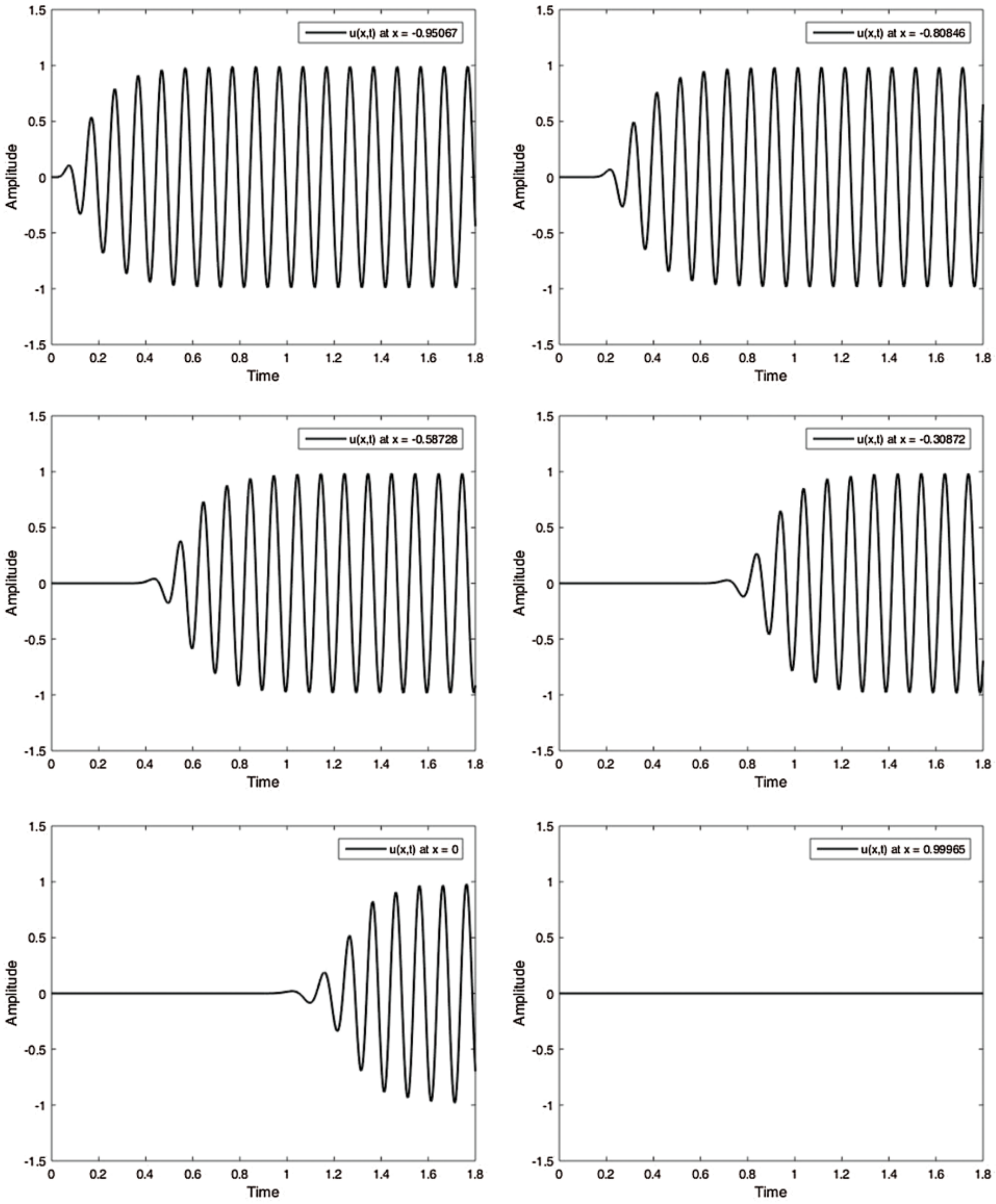

**Figure 6:** Eventual periodicity for linearized KdV equation: $u(x,t)$ at $x$=-0.950670, -0.808460, -0.587280, -0.308720, 0 and 0.999650 are shown through above six graphs respectively. For $\alpha$=1, $\beta$=0, $\mu$=$10^{-5}$, $x$∈[-1, 1], $C_{MQ}$=0.001, $N$=200, $N_x$=25, $tmax$=1.8, $\delta t$=0.001, $g(t)$=sin(20$\pi t$) tanh(5$t$)

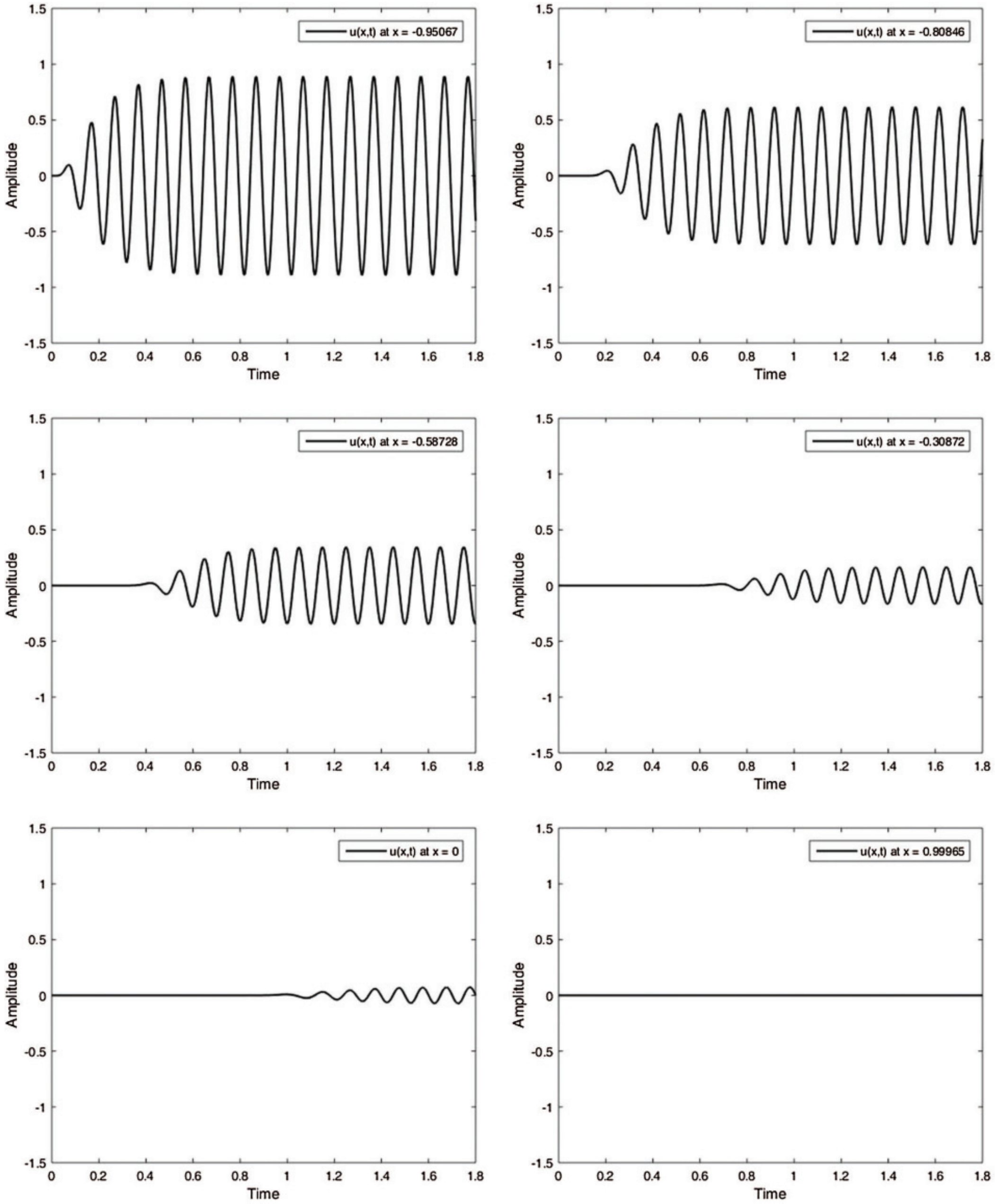

**Figure 7:** Eventual periodicity for linear KdV-Burgers equation: $u(x,t)$ at $x$=-0.950670, -0.808460, -0.587280, -0.308720, 0 and 0.999650 are shown through above six graphs respectively. For $\alpha$=1, $\beta$=0, $\mu$=$10^{-5}$ and $\delta$=$10^{-4}$, $x\in[-1, 1]$, $C_{MQ}$=0.001, $N$=200, $N_x$=25, $tmax$=1.8, $\delta t$=0.001, $g(t)$=sin(20$\pi$ $t$) tanh(5$t$)

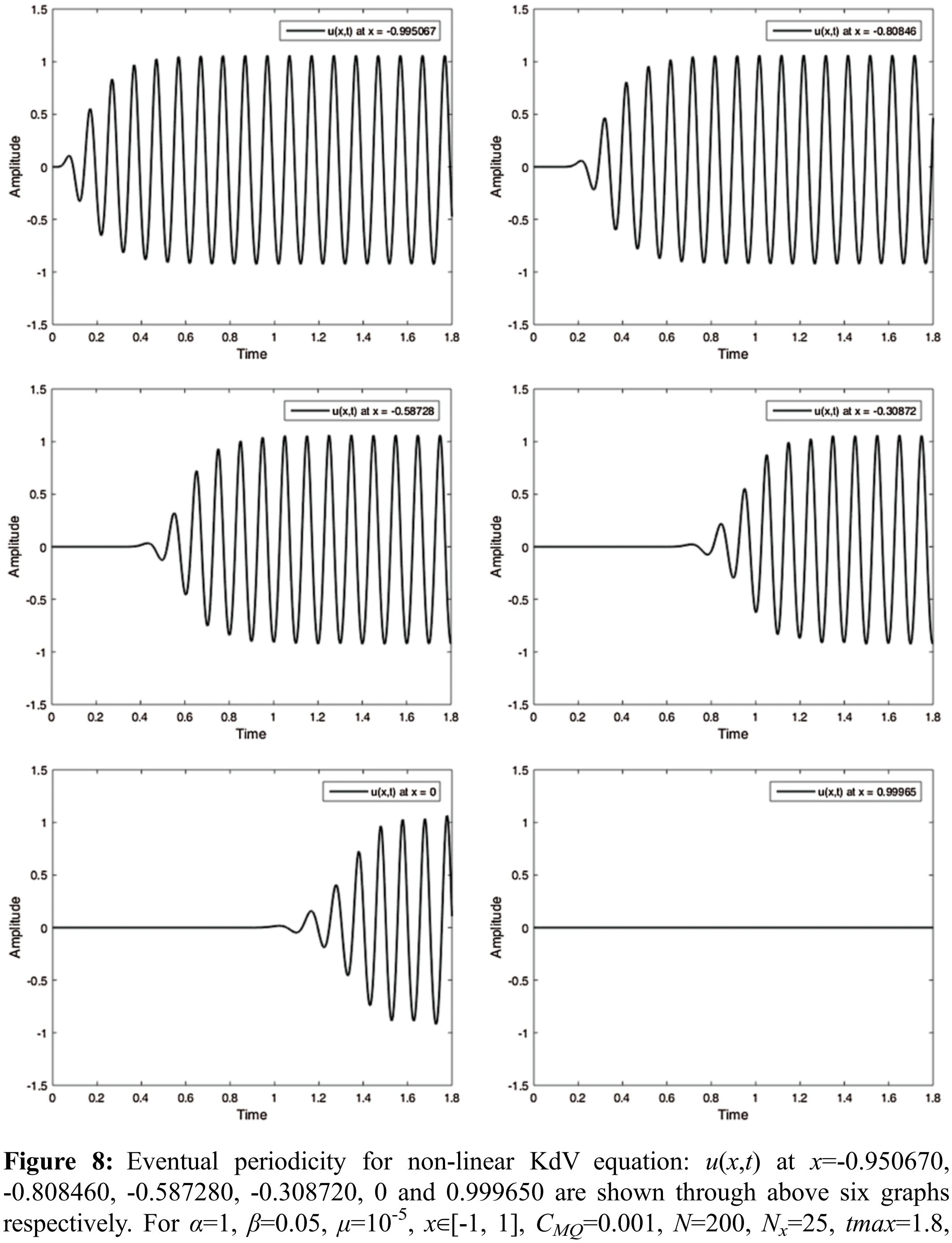

**Figure 8:** Eventual periodicity for non-linear KdV equation: $u(x,t)$ at $x$=-0.950670, -0.808460, -0.587280, -0.308720, 0 and 0.999650 are shown through above six graphs respectively. For $\alpha$=1, $\beta$=0.05, $\mu=10^{-5}$, $x\in[-1, 1]$, $C_{MQ}$=0.001, $N$=200, $N_x$=25, $tmax$=1.8, $\delta t$=0.001, $g(t)=\sin(20\pi t)\tanh(5t)$

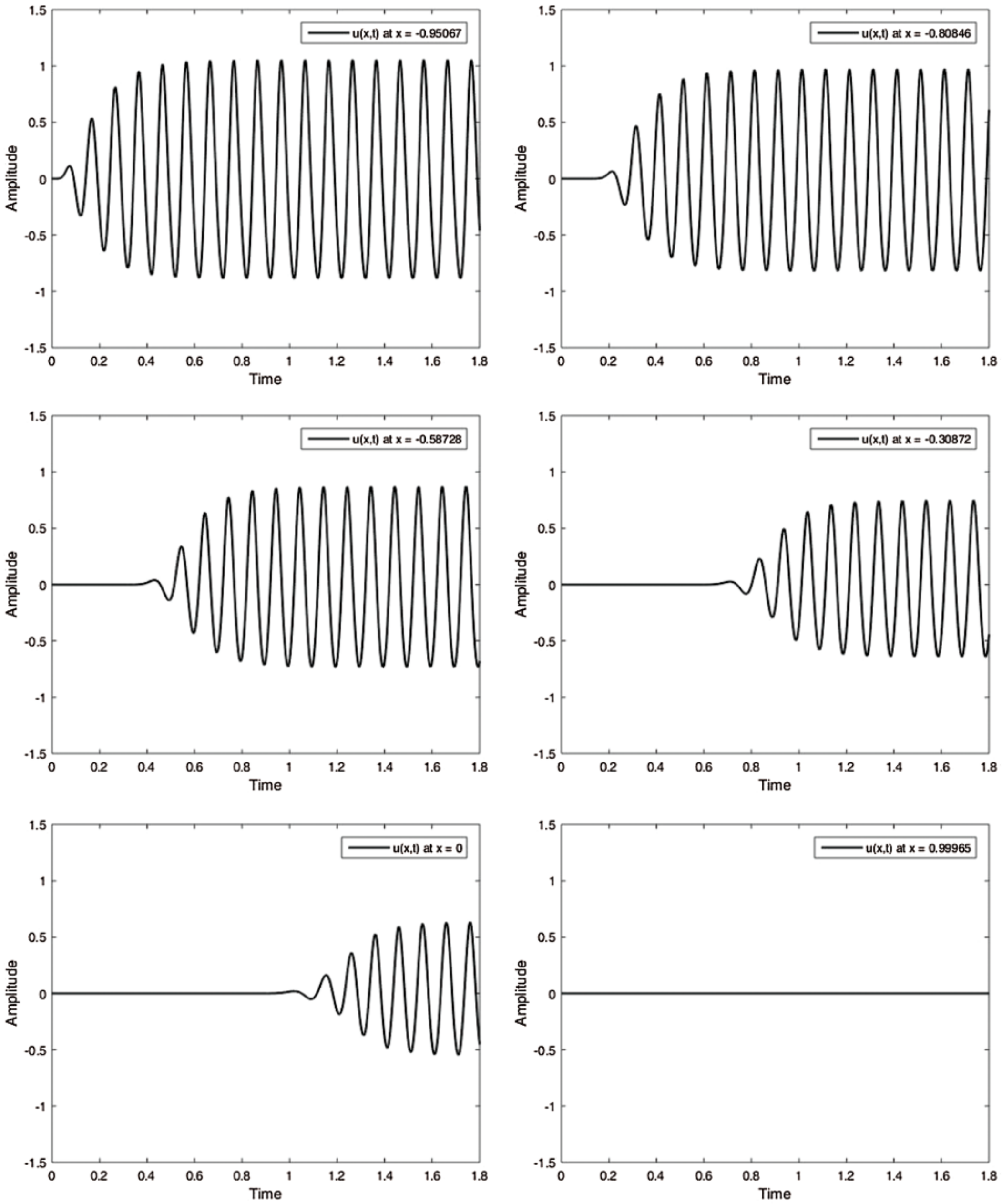

**Figure 9:** Eventual periodicity for KdV-Burgers equation: $u(x,t)$ at $x$=-0.950670, -0.808460, -0.587280, -0.308720, 0 and 0.999650 are shown through above six graphs respectively. For $\alpha$=1, $\beta$=0.05, $\mu$=$10^{-5}$, $\delta$=$10^{-4}$ $x\in[-1, 1]$, $C_{MQ}$=0.001, $N$=200, $N_x$=25, $tmax$=1.8, $\delta t$=0.001, $g(t)=\sin(20\pi t)\tanh(5t)$

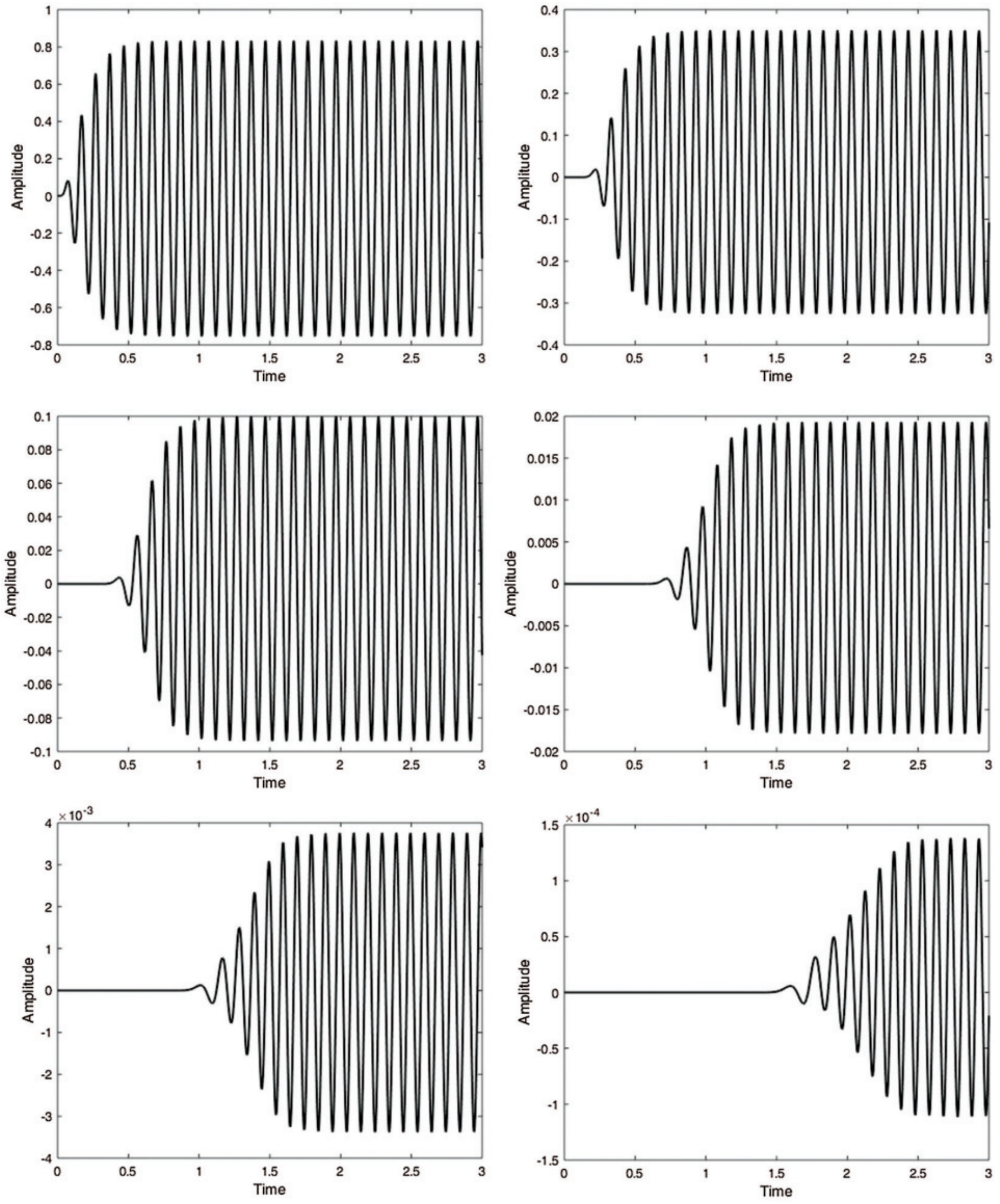

**Figure 10:** $u(x,t)$ at $x$=-0.950670, -0.808460, -0.587280, -0.308720, 0 and 0.999650 are shown through above six graphs respectively. For $\eta$=4.5, $\alpha$=1, $\beta$=0.05, $\mu=10^{-4}$, $x\in[-1, 1]$, $t\in[0, 3]$, $C_{MQ}$=0.001, $N$=200, $N_x$=25, $tmax$=3, $\delta t$=0.001, $g(t)=\sin(20\pi t)\tanh(5t)$

## 6 Concluding remarks

In the present work, the RBF-FD method is implemented to investigate numerically the eventual periodicity of solutions to the initial and boundary value problems (IBVPs) for some dispersive wave equations like Benjamin-Bona-Mahony and the Korteweg-de Vries on a bounded domain. We combined the localized RBF-FD numerical scheme with RK4 time integration scheme. The proposed numerical scheme is assessed for accuracy and efficiency by solving the KdV equation with a known analytic solution. Also the stability of the proposed numerical scheme is discussed and validated for KdV equation. The main benefit and important advantage of this RBF-FD scheme over other methods is its local and sparse nature of differentiation matrices, stability and higher order convergence rate and low computational cost. It is also applicable to other type of integral and non integral order PDEs.

**Acknowledgement:** The authors wish to express their appreciation to the reviewers for their helpful suggestions which greatly improved the presentation of this paper.

**Funding Statement**: The authors received no specific funding for this study.

**Conflicts of Interest**: The authors declare that they have no conflicts of interest to report regarding the present study.